\documentclass[MSNbibl,number,citesort,dvips]{arxbj}
\usepackage{upgreek}
\usepackage{graphicx}

\aid{0}
\volume{19}
\issue{5A}
\pubyear{2013}
\firstpage{1637}
\lastpage{1654}
\doi{10.3150/12-BEJ424} 

\makeatletter

\define@key{bid}{pmid}{}

\newtheorem{theorem}{Theorem}
\newtheorem{proposition}{Proposition}
\newtheorem{corollary}{Corollary}
\newremark{remark}{Remark}

\newcommand{\Faces}{\operatorname{Faces}}
\newcommand{\MaxFaces}{\operatorname{MaxPolytopes}}
\newcommand{\MaxFacets}{\operatorname{MaxPolytopes}_{d-1}}
\newcommand{\Cells}{\operatorname{Cells}}
\newcommand{\Vol}{\mathrm{V}}
\newcommand{\PHT}{\operatorname{PHT}}
\newcommand{\Cell}{\operatorname{Cell}}
\makeatother

\begin{document}
\begin{frontmatter}

\title{Geometry of iteration stable tessellations: Connection with
Poisson hyperplanes}
\runtitle{Geometry of STIT tessellations: Connection with Poisson hyperplanes}

\begin{aug}
\author[1]{\fnms{Tomasz} \snm{Schreiber}\thanksref{1,th1}}
\and
\author[2]{\fnms{Christoph} \snm{Th\"ale}\corref{}\thanksref{2}\ead[label=e2]{christoph.thaele@uni-osnabrueck.de}}
\runauthor{T. Schreiber and C. Th\"ale} 
\address[1]{Faculty of Mathematics and Computer Science, Nicolaus
Copernicus University, Toru\'n, Poland}
\address[2]{Institut f\"ur Mathematik, Universit\"at Osnabr\"uck,
Albrechtstra\ss e 28a, D-49076 Osnabr\"uck, Germany. \printead{e2}}
\end{aug}
\thankstext[*]{th1}{26.06.1975--01.12.2010}

\received{\smonth{4} \syear{2011}}
\revised{\smonth{12} \syear{2011}}

%
\begin{abstract}
Since the seminal work by Nagel and Weiss, the iteration stable (STIT)
tessellations have attracted considerable interest in stochastic
geometry as a natural and flexible, yet analytically tractable model
for hierarchical spatial cell-splitting and crack-formation processes.
We provide in this paper a fundamental link between typical
characteristics of STIT tessellations and those of suitable mixtures of
Poisson hyperplane tessellations using martingale techniques and
general theory of piecewise deterministic Markov processes (PDMPs). As
applications, new mean values and new distributional results for the
STIT model are obtained.
\end{abstract}

%

\end{frontmatter}

\section{Introduction}\label{introd}
Infinite divisibility or stochastic stability of a random object under
a certain operation is one of the most fundamental concepts in
probability theory. Prominent examples include the classical theory of
infinite divisible and stable distributions with their applications
around the central limit theorem, max-stable distributions studied in
extreme value theory or union infinitely divisible random sets studied
in classical stochastic geometry.

In the present paper, we deal with a class of iteration infinitely
divisible random tessellations of the $d$-dimensional Euclidean space
and, more specifically, with random tessellations that are \textit
{st}able under the operation of \textit{it}eration -- so-called STIT
tessellations. Recall that a tessellation (or mosaic) of ${\Bbb R}^d$
is a locally finite family of compact and convex polytopes with
pairwise no common interior points that cover the whole space. They are
one of the central objects studied in stochastic geometry and related
fields, see \cite{SW,SKM}. They are also of great importance for
applications of stochastic geometry to real-world problems for which we
refer to \cite{Beil0,Beil,Chi,Lau,Lau2}. In particular and as
discussed in~\cite{NMOW}, the STIT tessellations may serve as a
reference model for hierarchical spatial cell-splitting and crack
formation processes in natural sciences and technology, for example, to
describe geological or material phenomena or aging processes or surfaces.

The motivation for iteration stable tessellations can be traced back to
the 80s and the principle of iteration of tessellations can roughly be
explained as follows, cf. \cite{NW03,NW05}. Take a random primary or
frame tessellation and associate with each of its cells an independent
copy of the primary tessellation, called component tessellation, which
is independent of the primary tessellation as well. In each cell, a
local superposition of the primary tessellation and the associated
component tessellation is now performed. The described operation can be
applied repeatedly and we obtain this way a sequence of random
tessellations. It can be shown that, after appropriate rescaling, this
sequence converges to a limit tessellation, which is easily seen to be
infinitely divisible or even stable with respect to iteration
(depending on the stochastic properties of the primary tessellation).
%
\begin{figure}[b]

\includegraphics{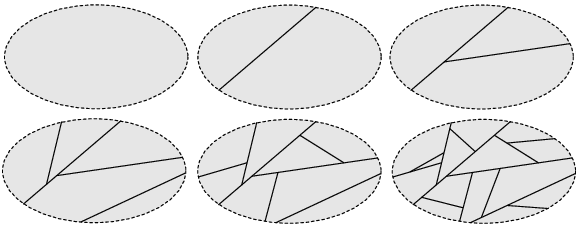}

\caption{Construction of a STIT tessellation in a convex window with
curved boundary.}\label{Fig0}
\end{figure}%

Starting with \cite{NW03}, STIT tessellations and their theoretical
framework were formally introduced in \cite{NW05}. In \cite{MNW,MNW2},
a tessellation-valued random process on the positive real
half-axis was constructed with the property that at each time the law
of the tessellation is stable under iteration. This dynamic point of
view also provides the link to a class of more general iteration
infinitely divisible tessellations, which are in the focus of the
present paper as well. In compact and convex windows $W\subset{\Bbb
R}^d$ with positive volume, this process can be explained as follows. A
terminal time $t>0$ and a (in some sense nondegenerate) measure
$\Lambda$ on the space of hyperplanes in ${\Bbb R}^d$ are fixed in
advance. Now $W$ is assigned a random lifetime. Upon expiry of its
lifetime, $W$ dies and splits into two random sets separated by a
hyperplane hitting $W$, which is chosen according to the suitable
normalization of $\Lambda$. The resulting new random sets are again
assigned independent random lifetimes and the entire construction
continues recursively until the deterministic time threshold $t$ is
reached, see Figure \ref{Fig0} for an illustration. The resulting
random structure tessellates the window $W$ and is denoted by
$Y(t\Lambda,W)$. In order to ensure the Markov property of this
construction in the continuous-time parameter $t$ and in order to keep
the law of $Y(t\Lambda,W)$ infinitely divisible or stable with respect
to iterations, we will have to take care of the special choice of the
lifetime distributions and the the cell-dividing hyperplanes, see
Section \ref{secSTIT} below. We would like to emphasize that the
dynamical representation is a special feature of STIT tessellations and
their infinitely divisible counterparts and, as recently pointed out in
\cite{Schr}, that such and similar spatio-temporal random processes
have remarkable potential for applications in stochastic geometry.

The purpose of the present paper is to explore the dynamic
representation further and to introduce a new technique, which unifies
and generalizes former approaches and which, moreover, has the
advantage that it allows to deal with properties of the model that were
out of reach so far. The crucial fact is that the construction above
has an interpretation as a piecewise deterministic Markov process
(PDMP) on the space of tessellations of $W$, which paves the way to the
general theory of PDMP's for which we refer to \cite{Jacobsen} in
particular. This in turn puts us into the position to construct certain
martingales related to the tessellation (see Section \ref
{secMARTINGALE}), eventually leading to fundamental comparison results
of the tessellations under consideration with certain mixtures of
Poisson hyperplane tessellations in Section \ref{secFIRST}. We would
like to remark at this point that the theory of PDMP's has previously
successfully been applied in stochastic geometry and spatial statistics
for the modeling of crack-growth networks, cf. \cite{Beil,Chi}. As
applications of our comparison results we calculate in Section \ref
{sectypifaces} several new mean values and also some new distributions
related to geometric objects determined by the tessellation. To keep
the paper self-contained, we recall in Section \ref{secLIMIT} the
construction of STIT tessellations as limit of repeated iterations and
formally introduce in Section \ref{secSTIT} their Markovian dynamic
representation. Moreover, we introduce there the above mentioned class
of iteration infinitely divisible random tessellations and summarize
some of their properties needed in this paper.

The current work is based on an extended version available online \cite{ST}.
It also forms the basis of our papers \cite{ST2,ST3,STP3,STP2}.
%
\section{STIT tessellations as limits}\label{secLIMIT}
Before explaining the concept of iteration of tessellations, let us fix
some basic notions and notation. A tessellation of ${\Bbb R}^d$ is a
locally finite partition of the space into compact convex polytopes,
the cells of the tessellation. One can regard a tessellation either as
a collection of its cells or as the closed set formed by the union of
their boundaries. We will follow the second point of view and denote by
$\Cells(Y)$ the set of cells of a tessellation $Y$. Thus, a random
tessellation can be regarded as a special random closed set in the
classical sense of stochastic geometry, see \cite{SW}. In particular,
this imposes the usual Fell topology and the corresponding Borel
measurable structure on the family of tessellations, see ibidem. By
$\mathcal{K}_d$, we denote in this paper the space of compact and
convex set in ${\Bbb R}^d$ with positive volume.

A random tessellation $Y$ is \textit{stationary} if its distribution
does not change upon actions of translations. Analogously a random
tessellation is said to be \textit{isotropic} if its distribution is
invariant under the action of the rotation group $\mathit{SO}_d$.

Whenever two random tessellations $Y_1$ and $Y_2$ of ${\Bbb R}^d$ are
given, we can define their \textit{iteration/nesting}. For this
purpose, we associate to each cell $c \in\Cells(Y_1)$ an independent
copy $Y_2(c)$ of $Y_2$ and we assume furthermore the family $\{Y_2(c)\dvt c\in\Cells(Y_1)\}$ to be independent of $Y_1$. Then we define
the iteration of $Y_1$ with $Y_2$ by
\[
Y_1\boxplus Y_2:=Y_1\cup\bigcup_{c\in\Cells(Y_1)}\bigl(Y_2(c)\cap c\bigr),
\]
that is, we take the local superposition of $Y_1$ and the family $\{Y_2(c)\dvt c\in\Cells(Y_1)\}$ inside the cells of $Y_1$. It was shown in
\cite{MNW2} that with $Y_1$ and $Y_2$ also $Y_1\boxplus Y_2$ is a
stationary random tessellation. A stationary random tessellation $Y$ is
called \textit{stable under iterations}, or \textit{STIT} for short, if
%
\begin{equation}\label{ITSTAB}
m( \underbrace{Y \boxplus\cdots\boxplus Y}_{m\ \mathrm{times}})\stackrel{D}{=} Y,\qquad m=2,3,\ldots,
\end{equation}
where $\stackrel{D}{=}$ stands for equality in distribution (note that
rescaling with factor $m$ ensures that the mean surface area of cell
boundaries per unit volume remains constant). In fact, using the
uniqueness results Theorem 3 and Corollary 2 in \cite{NW05} it is easy
to see that it is enough to take one fixed $m > 1$ in (\ref{ITSTAB}).

\begin{figure}

\includegraphics{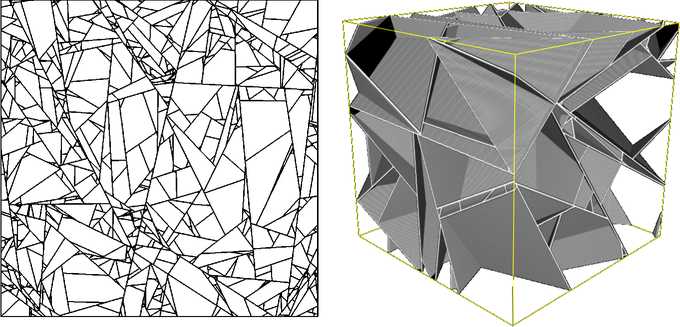}

\caption{Realizations of a planar and a spatial stationary and
isotropic STIT tessellation in a square and a cube,
respectively.}\label{Fig1}
\end{figure}

To proceed, let us be given a constant $0<t<\infty$ and an even
(symmetric) probability measure $\mathcal{R}$ on the unit sphere
${\Bbb S}^{d-1}$, usually identified with the induced distribution of
orthogonal hyperplanes on the space of $(d-1)$-dimensional linear
hyperplanes in ${\Bbb R}^d$, also denoted by $\mathcal{R}$ in the
sequel for notational simplicity. Define the measure $\Lambda$ on the
space $\mathcal{H}$ of affine hyperplanes in ${\Bbb R}^d$ by
%
\begin{equation}\label{LADEF}
\Lambda:=\ell_+ \otimes\mathcal{R},
\end{equation}
with $\ell_+$ standing for the Lebesgue measure on the positive real
half-axis $(0,\infty)$. Throughout this paper, we always require that
the support of $\mathcal{R}$ spans the whole space. Assume now that we
are given a stationary random tessellation $Y$ with surface intensity
$t$ (i.e., the mean surface area of cell boundaries per unit volume
equals $t$) and directional distribution $\mathcal{R}$ (i.e., the
distribution of the normal direction of the facet containing the
typical point is given by $\mathcal{R}$) and define the sequence
$(\mathcal{I}_n(Y))_{n\geq1}$ by
\[
\mathcal{I}_1(Y):=2(Y\boxplus Y), \qquad\mathcal{I}_{n}(Y):= \frac
{n}{n-1}\mathcal{I}_{n-1}(Y) \boxplus n Y = n (\underbrace{Y \boxplus
\cdots\boxplus Y}_{n \ \mathrm{times}}), \qquad n\geq2.
\]
It was shown in \cite{NW05}, Theorem 3, that $\mathcal{I}_n(Y)$ converges
in law, as $n \to\infty$, to a stationary random limit tessellation
$Y(t\Lambda)$ uniquely determined by $t\Lambda$. This tessellation is
easily shown to be stable under iterations, whence a STIT tessellation
with parameter $t$ and hyperplane measure $\Lambda$; we refer to
Figure \ref{Fig1} for an illustration of the limit tessellation.
%
%
\section{The Markovian construction}\label{secSTIT}
As already emphasized in the \hyperref[introd]{Introduction}, it is a crucial feature of
the random tessellations $Y(t\Lambda)$ introduced in the previous
section that they admit a simple and intuitive spatio-temporal
Markovian construction. For a restriction $Y(t\Lambda,W)$ of
$Y(t\Lambda)$ to a window $W\in\mathcal{K}_d$, this construction can
be described as follows (the reader is referred to \cite{NW05} for
full details). Assign to $W$ an exponentially distributed random
lifetime with parameter $\Lambda([W])$ where $[W] := \{H \in\mathcal
{H}\dvt H \cap W \neq\varnothing\}$ stands for the family of all
hyperplanes hitting $W$. Upon expiry of its lifetime, $W$ dies and
splits into $W^\pm=W\cap H^\pm$ ($H^\pm$ are the two half-spaces
determined by $H$), which are separated by a hyperplane in $[W]$ chosen
according to the law $\Lambda(\cdot\cap[W])/\Lambda([W])$. The
resulting new random sets $W^+$ and $W^-$ are again assigned
independent exponential lifetimes with respective parameters $\Lambda
([W^+])$ and $\Lambda([W^-])$ and the entire construction continues
recursively until the deterministic time threshold $t$ is reached, see
Figure \ref{Fig0}. The separating $(d-1)$-dimensional facets (the word
\textit{facet} stands for a $(d-1)$-dimensional face here and
throughout) arising in subsequent splits are usually referred to as
\textit{$(d-1)$-dimensional maximal polytopes} or \mbox{\textit{I}-segments} for $d=2$
as assuming shapes similar to the letter \textit{I}. (Note, that due
to a possibly curved boundary of $W$ some of the maximal polytopes may
also have a curved boundary and are no polytopes in the usual sense.
However, we abuse notation and include also these sets in our class of
maximal polytopes, which causes no difficulties in our theory.) The
resulting random closed set constructed inside $W$ is denoted by
$Y(t\Lambda,W)$, whereas the collection of all $(d-1)$-dimensional
maximal polytopes or \textit{I}-segments is denoted by $\MaxFacets(Y(t\Lambda
,W))$. Moreover, we write $\MaxFaces_k(Y(t\Lambda,W))$ for the
collection of \textit{$k$-dimensional maximal polytopes} of
$Y(t\Lambda,W)$, where a $k$-dimensional maximal polytope is just a
$k$-dimensional face of some $(d-1)$-dimensional maximal polytope
(again, some of them are no polytopes in the usual sense).

It was shown in \cite{NW05} that the law of $Y(t\Lambda,W)$ is
consistent in that $Y(t\Lambda,W) \cap V \stackrel{D}{=} Y(t\Lambda
,V)$ for convex $V \subset W$ and thus $Y(t\Lambda,W)$ can be extended
to a random tessellation $Y(t\Lambda)$ in the whole space, which is
then proved (see \cite{NW05}) to coincide with the limit tessellation
$Y(t\Lambda)$ considered in the previous section as notation already
suggests. Again, the family of all $k$-dimensional maximal polytopes of
$Y(t\Lambda)$ is denoted by $\MaxFaces_k(Y(t\Lambda))$ ($0\leq k\leq
d-1$). The stationary random tessellation $Y(t\Lambda)$ is
additionally isotropic if and only if $\mathcal{R}$ is the uniform
distribution on ${\Bbb S}^{d-1}$ in the factorization~(\ref{LADEF}).

A simple yet crucial observation is that even though only
translation-invariant measures $\Lambda$ of the form (\ref{LADEF})
show up in the limiting STIT tessellations, the dynamic construction
can be carried out with arbitrary non-atomic and locally finite
measures $\Lambda$ on $\mathcal{H}$ also leading to a consistent
family $Y(t\Lambda,W)$ and eventually, by extension, yielding a whole
space tessellation $Y(t\Lambda)$. Many of our theorems will be stated
in this general context. It should be emphasized that such
tessellations are no longer iteration stable (STIT). However, they have
the general property of being \textit{iteration infinitely divisible},
as they can be readily checked to arise as $m$-fold iterations of
$Y(t/m)$ for each $m \geq2$ in all finite windows. Formally, this
means that
\[
Y\bigl((t/m)\Lambda,W\bigr)^{\boxplus m} \stackrel{D}{=} Y(t\Lambda,W)
\]
for all $W\in\mathcal{K}_d$, which follows directly by construction
as yielding
\[
Y(s\Lambda,W) \boxplus Y(u\Lambda,W) \stackrel{D}{=} Y\bigl((s+u)\Lambda,W\bigr).
\]

It is worth pointing out that it is currently an open problem whether
\textit{any} iteration infinitely divisible random tessellation in
${\Bbb R}^d$ can be constructed in this way.
%
\section{Associated martingales}\label{secMARTINGALE}
The finite volume continuous-time incremental Markovian construction of
iteration infinitely divisible random tessellations, or more specially
of stationary STIT tessellations, as discussed in Sections \ref
{secLIMIT} and \ref{secSTIT} above, clearly enjoys the Markov property
in the continuous time parameter. Whence, natural martingales arise,
which will be of crucial importance for our further considerations. To
discuss these processes, we notice that for any $W\in\mathcal{K}_d$
and any nonatomic locally finite hyperplane measure $\Lambda$,
$(Y(t\Lambda,W))_{t>0}$ is a piecewise deterministic Markov process
(PDMP) on the space of tessellations of $W$ and we cite Chapter 3 in
\cite{GS}, Chapter 7 in \cite{Jacobsen} and Chapter 12 in
 \cite{Kallenberg} for the general theory of such processes. Using these
general results, we conclude that the PDMP $(Y(t\Lambda,W))_{t>0}$ has
its infinitesimal generator ${\Bbb L} := {\Bbb L}_{\Lambda; W}$ given by
%
\begin{equation}\label{GEN}
{\Bbb L}F(Y) = \int_{[W]} \sum_{f\in\Cells(Y\cap H)} [F(Y \cup\{f\}
) - F(Y)] \Lambda(\mathrm{d}H),
\end{equation}
where $Y$ is some instant of $Y(t\Lambda,W)$ and $F$ is some bounded
measurable function on space of tessellations of $W$, cf. Theorem 12.22
in \cite{Kallenberg} in particular. Notice that for a tessellation $Y$
of $W\in\mathcal{K}_d$ and $H\in[W]$, $\Cells(Y\cap H)$ stands for
the collection of $(d-1)$-dimensional cell of the sectional
tessellation $Y\cap H$. By standard theory as given in Lemma 5.1,
Appendix 1, Section 5 in \cite{KipLan} or, alternatively, by a direct
check we readily conclude now.
\begin{proposition}\label{propMART}
For $F$ bounded and measurable, the stochastic process
\[ \label{MART}\nonumber
F(Y(t\Lambda,W)) - \int_0^t {\Bbb L}F(Y(s\Lambda,W))\,\mathrm{d}s, \qquad t\geq0,
\]
is a martingale with respect to the filtration $\Im_t:=\sigma
(Y(s\Lambda,W)\dvt 0 \leq s \leq t)$.
\end{proposition}

To proceed towards the crucial Proposition \ref{propEXPECT}, consider
$F$ of the form
%
\begin{equation}\label{FDEF}
\Sigma_{\phi}(Y) := \sum_{f \in\MaxFacets(Y)} \phi(f),
\end{equation}
where $\phi(\cdot)$ is a generic bounded and measurable functional on
$(d-1)$-dimensional facets in $W$, that is to say a bounded and
measurable function on the space of closed $(d-1)$-dimensional
polytopes in $W$, possibly chopped off by the boundary of $W$, with the
standard measurable structure inherited from space of closed sets in $W$.
\begin{proposition}\label{propEXPECT} The stochastic process
\[
\Sigma_{\phi}(Y(t\Lambda,W)) - \int_0^t \int_{[W]} \sum_{f\in
\Cells(Y(s\Lambda,W)\cap H)} \phi(f) \Lambda(\mathrm{d}H)\,\mathrm{d}s,\qquad t\geq0,
\]
is a martingale with respect to $\Im_t$.
\end{proposition}
\begin{pf}
The functional $\Sigma_\phi$ defined by (\ref{FDEF}) is not
necessarily bounded and thus Proposition \ref{propMART} cannot be
applied directly with $F=\Sigma_\phi$ there. However, we can apply it
for the truncations $\Sigma_\phi^N := (\Sigma_\phi\wedge N) \vee
-N, N \in{\Bbb N}$ and let $N\to\infty$ to conclude that
\[
\Sigma_\phi(Y(t\Lambda,W))-\int_0^t {\Bbb L}\Sigma_\phi
(Y(s\Lambda,W))\,\mathrm{d}s,\qquad t\geq0
\]
is a \emph{local} $\Im_t$-martingale with a localizing sequence given
by
\[
\{\tau_N\}_{N\in{\Bbb N}} \qquad\mbox{with}\qquad \tau_N = \inf\{{t\geq
0}\dvt |\Sigma_\phi(Y(t\Lambda,W))| \geq N\}.
\]
Now, we apply the proof of Lemma 1 in \cite{NW05}, where the number of
cells in $Y(t\Lambda,W)$, and hence for all $Y(s\Lambda,W), s \leq
t$, has been shown to be bounded by a Furry--Yule-type linear birth
process whose cardinality at any given finite time admits moments of
all orders, to conclude that
\[
\Sigma_\phi(Y(t\Lambda,W))-\int_0^t {\Bbb L}\Sigma_\phi
(Y(s\Lambda,W))\,\mathrm{d}s, \qquad t \leq T,
\]
is of class DL for all $T>0$ in the sense of Definition 4.8 in
 \cite{KS}. Using now the result of Problem~5.19(i) ibidem, we finally
conclude that the random process
\[
\Sigma_\phi(Y(t\Lambda,W))-\int_0^t {\Bbb L}\Sigma_\phi
(Y(s\Lambda,W))\,\mathrm{d}s,\qquad t\geq0,
\]
is a martingale with respect to $\Im_t$. Moreover, in view of (\ref
{GEN}) we have
\begin{eqnarray*}
{\Bbb L}\Sigma_\phi(Y) &=& \int_{[W]}\sum_{f\in\Cells
(Y\cap H)}[\Sigma_\phi(Y\cup\{f\})- \Sigma_\phi(Y)] \Lambda(\mathrm{d}H)\\
&=& \int_{[W]}\sum_{f\in\Cells(Y\cap H)}\phi(f) \Lambda(\mathrm{d}H),
\end{eqnarray*}
which completes the proof.
\end{pf}
%
\section{Relationships for intensity measures}\label{secFIRST}
In this section, we establish two fundamental first-order properties of
$Y(t\Lambda,W)$ for general locally finite nonatomic measures $\Lambda
$, essentially obtained by comparison with suitable mixtures of Poisson
hyperplane tessellations. The results are somehow surprising as formal
identities and we are not able to provide an intuitive understanding.
However, their applications in Corollary \ref{corollary1} and in
Section \ref{sectypifaces} below have a very natural meaning and
interpretation.

The key to our results is Proposition \ref{propEXPECT} from the
previous section. To exploit it, consider the random measure
%
\begin{equation}\label{CellMeasure}
\mathcal{M}^{Y(t\Lambda,W)} := \sum_{c \in\Cells(Y(t\Lambda,W))}
\delta_c \quad\mbox{and}\quad {\Bbb M}^{Y(t\Lambda,W)} := {\Bbb E}\mathcal
{M}^{Y(t\Lambda,W)}
\end{equation}
with $\delta_c$ standing for the unit mass Dirac measure at $c$. In
full analogy, define $\mathcal{M}^{\PHT(t\Lambda,W)}$ and ${\Bbb
M}^{\PHT(t\Lambda,W)}$, where $\PHT(t\Lambda,W)$ is the Poisson
hyperplane tessellation with intensity measure $t\Lambda$, restricted
to $W$ (see \cite{SW} for background material). Further, put
\[\label{FK}
\mathcal{F}_k^{Y(t\Lambda,W)} := \sum_{f \in\MaxFaces_k(Y(t\Lambda
,W))} \delta_f,\qquad
{\Bbb F}_k^{Y(t\Lambda,W)} := {\Bbb E}\mathcal{F}_k^{Y(t\Lambda,W)},\quad
k=1,\ldots,d-1,
\]
where, recall, $\MaxFaces_k(Y)$ is the collection of $k$-dimensional
maximal polytopes of $Y$. Likewise, define
\[
\mathcal{F}_k^{\PHT(t\Lambda,W)} := \sum_{f \in\Faces_k(\PHT
(t\Lambda,W))}\!\delta_f, \qquad {\Bbb F}_k^{\PHT(t\Lambda,W)} := {\Bbb
E}\mathcal{F}_k^{\PHT(t\Lambda,W)},\quad k=1,\ldots,d-1,
\]
where $\Faces_k(\PHT(t\Lambda,W))$ is the collection of all $k$-face
of the Poisson hyperplane tessellation $\PHT(t\Lambda,W)$. Our first
claim is
\begin{theorem}\label{Meq} It holds that $ {\Bbb M}^{Y(t\Lambda,W)} =
{\Bbb M}^{\PHT(t\Lambda,W)} $.
\end{theorem}
\begin{pf}
Using (\ref{GEN}) and (\ref{MART}) with
\[
F(Y) := \sum_{c \in\Cells(Y)} \phi(c)
\]
for a general bounded measurable cell functional $\phi$, with a
localization argument as the one in the proof of Proposition \ref
{propEXPECT} we conclude that
%
\begin{eqnarray}\label{Mmart1}
&&\hspace*{-18pt}\int\phi(c') \mathcal{M}^{Y(t\Lambda,W)}(\mathrm{d}c') \nonumber\\
&&\hspace*{-18pt}\quad{}-
\int_0^t \int_{[W]} \sum_{f\in\Cells(Y(s\Lambda,W)\cap H)}
[\phi(\Cell^+(f,H|Y(s\Lambda,W)))+\phi(\Cell^-(f,H|Y(s\Lambda,W))) \\
&&\hspace*{-18pt}\hphantom{\quad{}-
\int_0^t \int_{[W]} \sum_{f\in\Cells(Y(s\Lambda,W)\cap H)}
[}{}- \phi(\Cell(f,H|Y(s\Lambda,W)))]
\Lambda(\mathrm{d}H) \,\mathrm{d}s,\qquad t\geq0,\nonumber
\end{eqnarray}
is a $\Im_t$-martingale. Here, $\Cell(f,H|Y(s\Lambda,W))$ stands for
the a.s. uniquely determined cell of $Y(s\Lambda,W)$ that gets divided
into $\Cell^{\pm}(f,H|Y(s\Lambda,W))$ by the facet $f$ on the
hyperplane $H\in[W]$. To simplify the notation, we use $c^{\pm}(H)$
to denote the cells into which $c$ gets divided by $H$, lying,
respectively, in the positive and negative half-space determined by
$H$. With this notation, (\ref{Mmart1}) says that
\begin{eqnarray*}
&&\int\phi(c') \mathcal{M}^{Y(t\Lambda,W)}(\mathrm{d}c')
\\
&&\quad{}- \int_0^t \int\int_{[c]} [\phi(c^+(H)) + \phi(c^-(H)) - \phi(c)]
\Lambda(\mathrm{d}H) \mathcal{M}^{Y(s\Lambda,W)}(\mathrm{d}c) \,\mathrm{d}s,\qquad t\geq0
\end{eqnarray*}
is a $\Im_t$-martingale. Taking expectations leads to
%
%
\begin{eqnarray}\label{EQNforM}
&& \int\phi(c') {\Bbb M}^{Y(t\Lambda,W)}(\mathrm{d}c')\nonumber
\\[-8pt]
\\[-8pt]
&&\quad = \int_0^t \int
\int_{[c]} [\phi(c^+(H)) + \phi(c^-(H)) - \phi(c)]
\Lambda(\mathrm{d}H) {\Bbb M}^{Y(s\Lambda,W)}(\mathrm{d}c) \,\mathrm{d}s\nonumber
\end{eqnarray}
for all bounded measurable $\phi$ as above.

To proceed, we regard ${\Bbb M}^{Y(s\Lambda,W)}$ as an element of the
space of bounded variation Borel
measures on the family of polyhedral sub-cells of $W$ endowed with the
standard measurable
structure inherited from the space of closed sets in $W$. Consider the linear
operator $T_{\Lambda}$ on this measure space given by
%
\begin{equation}\label{TM}
T_{\Lambda}\mu= \int\int_{[c]} [\delta_{c^+(H)} + \delta_{c^-(H)}
- \delta_{c}] \Lambda(\mathrm{d}H)
\mu(\mathrm{d}c).
\end{equation}
By (\ref{TM}), $ \|T_{\Lambda}\mu\|_{\mathrm{TV}} \leq(\int_{[W]}
\mathrm{d}\Lambda) \|\mu\|_{\mathrm{TV}}
= \Lambda([W]) \|\mu\|_{\mathrm{TV}}$ where $ \|\cdot\|_{\mathrm
{TV}}$ is the standard total variation norm of a measure. This
inequality turns into equality when $\mu= \delta_{W}$. Consequently,
$T_{\Lambda}$ is a bounded
operator of operator norm $\Lambda([W]) < +\infty$ by the assumed
locally finiteness of $\Lambda$.
Using the operator $T_\Lambda$, relation (\ref{EQNforM}) can be
rewritten in form of an initial value problem for the operator
differential equation
%
\begin{equation}\label{RNIEr}
\frac{\partial}{\partial t} {\Bbb M}^{Y(t\Lambda,W)} = T_{\Lambda}
{\Bbb M}^{Y(t\Lambda,W)},\qquad {\Bbb M}^{Y(0,W)} = \delta_{W},
\end{equation}
which, in view of the above properties of $T_{\Lambda}$, admits by
standard theory of linear operators (cf. \cite{Kato}, Chapter~IX.\S2, Section 2)
the unique
solution
%
\begin{equation}\label{ROZW}
{\Bbb M}^{Y(t\Lambda,W)} = \exp(t T_{\Lambda}) \delta_{W},\qquad t \geq0,
\end{equation}
where the operator exponential of $tT_\lambda$ is applied to the
measure $\delta_W$. It is easily seen that exactly the same equations
(\ref{EQNforM}), (\ref{RNIEr}) and thus also
(\ref{ROZW}) hold for ${\Bbb M}^{\PHT(t\Lambda,W)}$. In particular,
we have ${\Bbb M}^{Y(t\Lambda,W)} = {\Bbb M}^{\PHT(t\Lambda,W)}$, as
required.
\end{pf}

It is interesting to note that in the translation-invariant set-up, we
obtain as a corollary that the distribution ${\Bbb Q}$ of the typical
cell of the STIT tessellation $Y(t\Lambda)$ coincides with the typical
cell distribution ${\Bbb Q}^{\PHT(t\Lambda)}$ of a stationary Poisson
hyperplane tessellation $\PHT(t\Lambda)$ with intensity measure
$t\Lambda$, a result that has previously been shown in \cite{NW03} by
completely different arguments. Recall, that -- in intuitive terms --
the typical cell of a tessellation is a randomly selected cell, where
each cell has the same chance of being selected.
\begin{corollary}\label{corollary1}
In the translation-invariant set-up as considered in the previous
paragraph it holds that $\Bbb Q={\Bbb Q}^{\PHT(t\Lambda)}$.
\end{corollary}
\begin{pf}
The distribution ${\Bbb Q}$ is formally defined by the relation
%
\begin{equation}\label{eqdistqtypcell}
\lambda_d{\Bbb Q}(A) = \lim_{r\rightarrow\infty}{{\Bbb E}\sum
_{c\in\Cells(Y(t\Lambda))}\boldsymbol{1}[c\subset rW]\boldsymbol
{1}[c-m(c)\in A]\over r^d\Vol_d(W)},
\end{equation}
where $A$ is a measurable subset in the space of $d$-dimensional
polytopes, $\lambda_d$ is the cell density of $Y(t\Lambda)$, $W\in
\mathcal{K}_d$ and where $m(c)$ stands for some translation-covariant
selector of the $d$-dimensional polytope $c$ (for example the Steiner
point or the center of gravity), cf. \cite{SW}, equation (4.8, 4.9). The
distribution ${\Bbb Q}^{\PHT(t\Lambda)}$ is defined in a similar
spirit. Here and below $\boldsymbol{1}[\cdot]$ stands for the usual
indicator function, which is $1$ if the statement in brackets if
fulfilled and $0$ otherwise. Rewriting the sum in (\ref
{eqdistqtypcell}) as an integral and using Campbell's theorem, we obtain
\begin{eqnarray*}
&&{\Bbb E}\sum_{c\in\Cells(Y(t\Lambda))}\boldsymbol{1}[c\subset
rW]\boldsymbol{1}[c-m(c)\in A]
\\
&&\quad=\int\boldsymbol{1}[c\subset rW]\boldsymbol{1}[c-m(c)\in A]{\Bbb
M}^{Y(t\Lambda,rW)}(\mathrm{d}c),
\end{eqnarray*}
where the condition $c\subset rW$ excludes from $Y(t\Lambda
,rW)\stackrel{D}{=}Y(t\Lambda)\cap rW$ those cells that hit the
boundary of $rW$. Theorem \ref{Meq} allows now to replace ${\Bbb
M}^{Y(t\Lambda,rW)}$ by ${\Bbb M}^{\PHT(t\Lambda,rW)}$. Moreover,
Theorem \ref{Meq} clearly implies that the tessellations $Y(t\Lambda
)$ and $\PHT(t\Lambda)$ have the same cell density, which in view of
(\ref{eqdistqtypcell}) completes the argument.
\end{pf}

Having characterized ${\Bbb M}^{Y(t\Lambda,W)}$, we now turn to the
lower-dimensional face intensity measures ${\Bbb F}_k^{Y(t\Lambda,W)}$
in our general set-up, not necessarily assuming the
translation-invariance of the hyperplane measure $\Lambda$.
\begin{theorem}\label{FKeq}
For all $k=0,\ldots,d-1$, it holds that
\[
{\Bbb F}_k^{Y(t\Lambda,W)} = (d-k) 2^{d-k-1} \int_0^t \frac{1}{s}
{\Bbb F}_k^{\PHT(s\Lambda,W)} \,\mathrm{d}s.
\]
\end{theorem}
\begin{pf}
Fix $k \in\{0,\ldots,d-1\}$. Let $\psi$ be a general bounded
measurable function of a
$k$-dimensional maximal polytope, as usual regarded as a closed subset
of $W$,
and for a $(d-1)$-dimensional maximal polytope $h$ put
%
\begin{equation}\label{PhiPsi}
\phi(h) := \sum_{f \in\Faces_k(h)} \psi(f) ,
\end{equation}
noting that the $k$-dimensional maximal polytopes of the tessellation
$Y(t\Lambda,W)$ are precisely the $k$-faces of its $(d-1)$-dimensional
maximal polytopes.
Using Proposition \ref{propEXPECT}, taking expectations and recalling
(\ref{CellMeasure}) we see that
\[
{\Bbb E} \Sigma_{\phi}(Y(t\Lambda,W)) = \int\phi \,\mathrm{d}{\Bbb
F}_{d-1}^{Y(t\Lambda,W)} =
\int_0^t \int\int_{[c]} \phi(c \cap H) \Lambda(\mathrm{d}H) {\Bbb
M}^{Y(s\Lambda,W)}(\mathrm{d}c) \,\mathrm{d}s.
\]
Applying Theorem \ref{Meq}, we get
%
\begin{equation}\label{ConclAux1}
\int\phi \,\mathrm{d}{\Bbb F}_{d-1}^{Y(t\Lambda,W)} = \int_0^t \int\int
_{[c]} \phi(c \cap H)
\Lambda(\mathrm{d}H) {\Bbb M}^{\PHT(s\Lambda,W)}(\mathrm{d}c) \,\mathrm{d}s.
\end{equation}
Now the Slivnyak--Mecke formula \cite{SW}, Theorem 3.2.5 implies
%
\begin{equation}\label{eqhinzugefuegt}
\int\int_{[c]} \phi(c \cap H) \Lambda(\mathrm{d}H) {\Bbb M}^{\PHT(\Lambda,W)}(\mathrm{d}c)
= \int\phi \,\mathrm{d}{\Bbb F}_{d-1}^{\PHT(\Lambda,W)}.
\end{equation}
Indeed, identifying $\PHT(\Lambda,W)$ with the collection of Poisson
hyperplanes hitting $W$ we have
\begin{eqnarray*}
&&\int\phi \,\mathrm{d}{\Bbb F}_{d-1}^{\PHT(\Lambda,W)}\\
&&\quad= {\Bbb
E}\sum_{H\in\PHT(\Lambda,W)}\sum_{f\in\Cells(\PHT(\Lambda
,W)\cap H)}\phi(f)= \int_{[W]}{\Bbb E}\sum_{f\in\Cells(\PHT(\Lambda
,W)\cap H)}\phi(f) \Lambda(\mathrm{d}H)
\\
&&\quad={\Bbb E}\sum_{c\in\Cells(\PHT(\Lambda,W))}\int
_{[c]}\phi(c\cap H) \Lambda(\mathrm{d}H)=\int\int_{[c]} \phi
(c \cap H) \Lambda(\mathrm{d}H) {\Bbb M}^{\PHT(\Lambda,W)}(\mathrm{d}c).
\end{eqnarray*}
Taking in (\ref{eqhinzugefuegt}) now $s\Lambda$ in place of $\Lambda
$ we find more generally
\[
\int\int_{[c]} \phi(c \cap H) \Lambda(\mathrm{d}H) {\Bbb M}^{\PHT(s\Lambda,W)}(\mathrm{d}c)
= \frac{1}{s} \int\phi \,\mathrm{d}{\Bbb F}_{d-1}^{\PHT(s\Lambda,W)},
\]
whence, with (\ref{ConclAux1}),
\[
\int\phi \,\mathrm{d}{\Bbb F}_{d-1}^{Y(t\Lambda,W)} = \int_0^t \frac{1}{s}
\int\phi \,\mathrm{d}{\Bbb F}_{d-1}^{\PHT(s\Lambda,W)} \,\mathrm{d}s.
\]
We note now that, by (\ref{PhiPsi}),
\[
\int\phi \,\mathrm{d}{\Bbb F}_{d-1}^{Y(t\Lambda,W)} = \int\psi \,\mathrm{d}{\Bbb F}_k^{Y(t\Lambda,W)},
\]
because each $k$-dimensional maximal polytope is a $k$-face of
precisely one $(d-1)$-dimensional maximal polytope in $Y(t\Lambda,W)$.
On the other hand,
\[
\int\phi \,\mathrm{d}{\Bbb F}_{d-1}^{\PHT(s\Lambda,W)} = (d-k) 2^{d-k-1} \int
\psi \,\mathrm{d}{\Bbb F}_k^{\PHT(s\Lambda,W)},
\]
because each $k$-face of $\PHT(s\Lambda,W)$ is a $k$-face of $(d-k)
2^{d-k-1}$ facets of $\PHT(s\Lambda,W)$, see Theorems 10.1.2 and
10.3.1 in \cite{SW}.
Hence, we conclude that
\[
\int\psi \,\mathrm{d}{\Bbb F}_{k}^{Y(t\Lambda,W)} = (d-k) 2^{d-k-1} \int_0^t
\frac{1}{s} \int\psi \,\mathrm{d}{\Bbb F}_k^{\PHT(s\Lambda,W)} \,\mathrm{d}s
\]
for all $\psi$ bounded and measurable, which completes the proof of
the theorem.
\end{pf}

Some of our arguments in the sequel and the theory developed in
 \cite{ST2,ST3}, for example, require a straightforward formal extension of
Theorem \ref{FKeq}. Namely, we formally \textit{mark} all
$(d-1)$-dimensional maximal polytopes of the tessellation $Y(t\Lambda
,W)$ by their \textit{birth times}. This gives rise to the birth-time
augmented tessellation $\hat{Y}(t,W)$ with birth-time-marked
$(d-1)$-dimensional maximal polytopes and makes the Markovian
construction of $\hat{Y}(t,W)$ into a Markov process whose generator
$\hat{{\Bbb L}}$ is a clear modification of ${\Bbb L}$ as given in
(\ref{GEN}):
\[
\hat{\Bbb L} \hat{F}(\hat{Y}(s,W)) = \int_{[W]} \sum_{f\in\Cells
(Y(s,W)\cap H)} \bigl[\hat{F}(\hat{Y}(s,W) \cup[\{f\},s]) - \hat{F}(\hat
{Y}(s,W))\bigr] \Lambda(\mathrm{d}H)
\]
for $\hat{F}$ bounded and measurable on the space of birth-time-marked
tessellations of $W$. Consequently, writing $\hat{\Bbb F}^{Y(t\Lambda
,W)}_k, k=0,\ldots,d-1$, for the birth-time-marked version of
${\Bbb F}^{Y(t\Lambda,W)}_k$, where each $k$-dimensional maximal
polytope is marked with its birth time, by
a straightforward modification of the proof of Theorem \ref{FKeq} we
are led to
the following.
\begin{corollary}\label{TIMEMARKEDFQeq}
For all $k=0,\ldots,d-1$ it holds that
\[
\hat{\Bbb F}^{Y(t\Lambda,W)}_k = (d-k) 2^{d-k-1} \int_0^t \frac{1}{s}
\bigl[ {\Bbb F}^{\PHT(s\Lambda,W)}_k \otimes\delta_s \bigr] \,\mathrm{d}s.
\]
\end{corollary}
%
\section{Typical maximal polytope distributions}\label{sectypifaces}
We are now going to apply the results obtained in the last section to
the stationary set-up, that is, with $\Lambda$ translation-invariant,
to study the distribution of typical $k$-dimensional maximal polytopes
of the STIT tessellation $Y(t\Lambda)$ in ${\Bbb R}^d$. Recall, that,
in intuitive terms, the typical $k$-dimensional maximal polytope of
$Y(t\Lambda)$ is what we get when we equiprobably choose one of the
tessellations $k$-dimensional maximal polytopes. The typical
$k$-dimensional maximal polytope distribution ${\Bbb Q}_k$ ($1\leq
k\leq d-1$) is formally given by
%
\begin{equation}\label{eqtypff}
\lambda_k{\Bbb Q}_k(A)=\lim_{r\rightarrow\infty}{1\over\Vol
_d(rW)}\sum_{f\in\MaxFaces_k(Y(t\Lambda))}\boldsymbol{1}[f\subset
rW]\boldsymbol{1}[f-m(f)\in A],
\end{equation}
where $\lambda_k$ is the mean number of $k$-dimensional maximal
polytope selectors $m(f)$ of $Y(t\Lambda)$ per unit volume (this is
the intensity of $\MaxFaces_k(Y(t\Lambda))$), $A$ is a measurable
subset of the space of $k$-dimensional polytopes, $W\in\mathcal{K}_d$
and, as in the proof of Corollary \ref{corollary1}, where $m(f)$
stands for some translation-covariant selector of $f$ (as its Steiner
point for example). Similarly, the distribution ${\Bbb Q}_k^{\PHT
(t\Lambda)}$ of the typical $k$-face of a Poisson hyperplane
tessellation with intensity measure $t\Lambda$ is defined, see Chapter
4.1 in \cite{SW}.
\begin{theorem}\label{CORDISTR} For $k\in\{1,\ldots,d-1\}$, the
distribution ${\Bbb Q}_k$ of the typical $k$-dimensional maximal
polytope of $Y(t\Lambda)$ is given by
\[
{\Bbb Q}_k=\int_0^t{{d}s^{d-1}\over t^d}{\Bbb Q}_k^{\PHT(s\Lambda)}\,\mathrm{d}s.
\]
\end{theorem}

The preceding theorem can be rephrased by saying that the distribution
of the typical $k$-dimensional maximal polytope of a STIT tessellation
is a mixture of suitable rescalings of distributions of typical
$k$-dimensional faces of Poisson hyperplane tessellations, and that the
mixing distribution is a beta-distribution on $(0,t)$ with parameters
$d$ and $1$, which has density ${\,\mathrm{d}s^{d-1}\over t^d}\boldsymbol{1}[0<s<t]$.
\begin{pf*}{Proof of Theorem \ref{CORDISTR}}
To start, let $\varphi_k$ be a real-valued, bounded,
translation-invariant, non-negative measurable function on the space of
$k$-dimensional polytopes ($1\leq k\leq d-1$) and denote by $\overline
{\varphi}_k(Y(t\Lambda))$ the (possibly infinite) $\varphi
_k$-density of $Y(t\Lambda)$ in the sense of \cite{SW}, Chapter 4.1,
that is,
%
\begin{equation}\label{densitydefqq}
\overline{\varphi}_k(Y(t\Lambda))=\lim_{r\rightarrow\infty
}{1\over\Vol_d(rW)}{\Bbb E}\sum_{f\in\MaxFaces_k(Y(t\Lambda
))}\boldsymbol{1}[f\subset rW]\varphi_k(f)
\end{equation}
with $W\in\mathcal{K}_d$. The existence of this limit is guaranteed
by Theorem 4.1.3 ibidem. Using now Campbell's theorem \cite{SW}, Theorem~3.1.2 and Theorem \ref{FKeq} from above, we obtain, possibly with
both sides infinite,
\begingroup
\abovedisplayskip=6.7pt
\belowdisplayskip=6.7pt
\begin{eqnarray*}
\overline{\varphi}_k(Y(t\Lambda))
&=& \lim_{r\rightarrow
\infty}{1\over\Vol_d(rW)}{\Bbb E}\sum_{f\in\MaxFaces_k(Y(t\Lambda
))}\boldsymbol{1}[f\subset rW]\varphi_k(f)\\[-1pt]
&=& \lim_{r\rightarrow\infty}{1\over\Vol_d(rW)}{\Bbb
E}\int\boldsymbol{1}[f\subset rW]\varphi_k(f) \mathcal
{F}_k^{Y(t\Lambda,rW)}(\mathrm{d}f)\\[-1pt]
&=& \lim_{r\rightarrow\infty}{1\over\Vol_d(rW)}\int
\boldsymbol{1}[f\subset rW]\varphi_k(f) {\Bbb F}_k^{Y(t\Lambda
,rW)}(\mathrm{d}f)\\[-1pt]
&=& \lim_{r\rightarrow\infty}{(d-k)2^{d-k-1}\over\Vol
_d(rW)}\int_0^t{1\over s}\int\boldsymbol{1}[f\subset rW]\varphi
_k(f) {\Bbb F}_k^{\PHT(s\Lambda,rW)}(\mathrm{d}f)\,\mathrm{d}s.
\end{eqnarray*}
By dominated convergence, we can continue as follows:
%
\begin{eqnarray}\label{EQMIXFK}
&=& (d-k)2^{d-k-1}\int_0^t{1\over s} \biggl[\lim_{r\rightarrow
\infty}{1\over\Vol_d(rW)}{\Bbb E}\sum_{f\in\Faces_k(\PHT
(s\Lambda))}\boldsymbol{1}[f\subset rW]\varphi_k(f) \biggr]\,\mathrm{d}s\nonumber
\\[-9pt]
\\[-9pt]
&=& (d-k)2^{d-k-1}\int_0^t{1\over s}\overline{\varphi}_k(\PHT
(s\Lambda)) \,\mathrm{d}s.\nonumber
\end{eqnarray}
For $\varphi_k\equiv1$ we find $\lambda_k=\overline{\varphi
}_k(Y(t\Lambda))$, which clearly has an interpretation as intensity of
$k$-dimensional maximal polytopes. Similarly, we denote by $\lambda
_k^{\PHT(s\Lambda)}$ the intensity of $k$-faces of the Poisson
hyperplane tessellation $\PHT(s\Lambda)$ and notice that $s\PHT
(s\Lambda)\stackrel{D}{=}\PHT(\Lambda)$. Hence, there is a constant
$C_k\in(0,\infty)$, which is independent of $s$, such that $\lambda
_k^{\PHT(s\Lambda)}=C_ks^d$ (in fact $C_k$ is explicitly known and
given by equation (10.44) in \cite{SW}). Thus, in view of (\ref{EQMIXFK})
we get
%
\begin{equation}\label{NKIEXPR}
\lambda_k=(d-k)2^{d-k-1}\int_0^t{1\over s}C_ks^d\,\mathrm{d}s={d-k\over
d}2^{d-k-1}C_kt^d={d-k\over d}2^{d-k-1}\lambda_k^{\PHT(t\Lambda)}.
\end{equation}
Let $m(f)$ with $f\in\MaxFaces_k(Y(t\Lambda))$ be as in (\ref
{eqtypff}). Using (\ref{EQMIXFK}), this time with $\varphi
_k(f):=\boldsymbol{1}[f-m(f)\in\cdot]$, together with (\ref
{eqtypff}), we obtain
%
\begin{equation}\label{MIXTYPEQ}
\lambda_k{\Bbb Q}_k=(d-k)2^{d-k-1}\int_0^t{1\over s}\lambda_k^{\PHT
(s\Lambda)}{\Bbb Q}_k^{\PHT(s\Lambda)}\,\mathrm{d}s.
\end{equation}
Combining (\ref{MIXTYPEQ}) with (\ref{NKIEXPR}) we find, with $C_k$
as above,
\begin{eqnarray*}
{\Bbb Q}_k &=& (d-k)2^{d-k-1}\int_0^t{1\over s}{\lambda
_k^{\PHT(s\Lambda)}\over\lambda_k}{\Bbb Q}_k^{\PHT(s\Lambda)}\,\mathrm{d}s =
\int_0^t{d\over s}{\lambda_k^{\PHT(s\Lambda)}\over\lambda_k^{\PHT
(t\Lambda)}}{\Bbb Q}_k^{\PHT(s\Lambda)}\,\mathrm{d}s\\[-1pt]
&=& \int_0^t{d\over s}{C_ks^d\over C_kt^d}{\Bbb Q}_k^{\PHT
(s\Lambda)}\,\mathrm{d}s = \int_0^t{{d}s^{d-1}\over t^d}{\Bbb Q}_k^{\PHT(s\Lambda)}\,\mathrm{d}s,
\end{eqnarray*}
which is the desired expression for ${\Bbb Q}_k$.
\end{pf*}\eject
\endgroup

It is interesting to note that formally marking the $k$-dimensional
maximal polytopes with their birth-times and repeating the argument
leading to Theorem \ref{CORDISTR} with Theorem \ref{FKeq} replaced by
its time-marked extension in Corollary \ref{TIMEMARKEDFQeq} we obtain
the birth time-marked extension of Theorem \ref{CORDISTR}:
\begin{corollary}\label{TIMEMARKEDCORDISTR}
The distribution $\hat{\Bbb Q}_k$ of the typical birth-time-marked
$k$-dimensional maximal polytope of $Y(t\Lambda)$ is given by
\begingroup
\abovedisplayskip=6.7pt
\belowdisplayskip=6.7pt
\[
\hat{\Bbb Q}_k=\int_0^t{ds^{d-1}\over t^d}\bigl [ {\Bbb Q}_k^{\PHT
(s\Lambda)} \otimes\delta_s \bigr] \,\mathrm{d}s,
\]
where $k\in\{1,\ldots,d-1\}$ as above.
\end{corollary}

From these identities, mean values and also distributional results for
the typical $k$-dimensional maximal polytope of a stationary STIT
tessellation can be deduced. To illustrate the general method, we
exemplarily calculate at first the mean intrinsic volumes of the
typical $k$-dimensional maximal polytope of the STIT tessellation
$Y(t\Lambda)$. To neatly formulate them, let $\Pi$ be the associated
zonoid of the Poisson hyperplane tessellation $\PHT(t\Lambda)$. This
is the convex body which has its support $h_\Pi(u)$ function given by
\[
h_\Pi(u)={t\over2}\int_{{\Bbb S}^{d-1}}|\langle u,v\rangle|\mathcal
{R}(\mathrm{d}v), \qquad u\in{\Bbb R}^d,
\]
see \cite{SW}, equation (4.59), where $\langle\cdot,\cdot\rangle$
denotes the standard scalar product in ${\Bbb R}^d$. Moreover, we will
denote by $V_j(K)$ the intrinsic volume of order $j$ of $K\in\mathcal
{K}_d$ in the usual sense of integral geometry. In particular, $V_d(K)$
is the volume, $2V_{d-1}(K)$ the surface area, $V_1(K)$ a constant
multiple of the mean width and $V_0(K)=1$ the Euler-characteristic of
$K$, cf.~\cite{SW}.
\begin{corollary}\label{corintvols}
Fix $0\leq j\leq k\leq d-1$. The mean $j$th intrinsic volume $V_j$ of
the typical $k$-dimensional maximal polytope $I_k$ (this is a random
polytope with distribution ${\Bbb Q}_k$) is given by
\[
{\Bbb E}V_j(I_k)={d\over d-j}{{d-j\choose d-k}\over{d\choose
k}}{V_{d-k}(\Pi)\over\Vol_d(\Pi)}.
\]
\end{corollary}
\begin{pf}
Using Theorem \ref{CORDISTR}, Theorem 10.3.3 in \cite{SW} and the
homogeneity of the intrinsic volumes, we get
\begin{eqnarray*}
{\Bbb E}V_j(I_k)
&=& \int_0^t{ds^{d-1}\over t^d}{{d-j\choose
d-k}V_{d-j} ((s/t)\Pi)\over{d\choose k}\Vol_d ((s/t)\Pi
)}\,\mathrm{d}s={{d-j\choose d-k}\over{d\choose k}}{V_{d-k}(\Pi)\over\Vol
_d(\Pi)}\int_0^t{ds^{d-1}\over t^d}{(s/t)^{d-j}\over
(s/t)^d}\,\mathrm{d}s\\
&=& {d\over d-j}{{d-j\choose d-k}\over{d\choose
k}}{V_{d-k}(\Pi)\over\Vol_d(\Pi)},
\end{eqnarray*}
which completes the proof.
\end{pf}\eject
\endgroup

Note that in the isotropic case, that is, when $\mathcal{R}$ is the
uniform distribution on ${\Bbb S}^{d-1}$ in the factorization (\ref
{LADEF}), the zonoid $\Pi$ is a $d$-dimensional ball with radius
proportional to $t$. More specifically, we have in this case
\[
{\Bbb E}V_j(I_k)={d\over d-j}{k\choose j} \biggl({d\kappa_d\over\kappa
_{d-1}} \biggr)^j{1\over\kappa_jt^j},
\]
where $\kappa_j$ is the volume of the $j$-dimensional unit ball.
\begin{remark}
In the planar ($d=2$) and in the spatial case ($d=3$), the mean values
${\Bbb E}V_j(I_k)$ are in accordance with the values obtained earlier
in \cite{MNW07,NW08}, for example. The method there was based on the
stochastic stability of the tessellation under iterations, which leads
to balance equations for ${\Bbb E}V_j(I_k)$ that can be solved by using
intersection formulae for random tessellations. It seems, however, that
this method becomes impracticable in higher space dimensions.
\end{remark}
\begin{remark}
Corollary \ref{corintvols} shows that the intrinsic volumes $V_j$
($0\leq j\leq k$) are integrable with respect to ${\Bbb Q}_k$ ($1\leq
k\leq d-1$), the typical $k$-dimensional maximal polytope distribution.
In view of Theorem 4.1.2 in \cite{SW} this allows us to replace the
definition (\ref{densitydefqq}) of $\overline{\varphi}_k$ by
\begin{eqnarray*}
\overline{\varphi}_k(Y(t\Lambda))
&=& \lim_{r\rightarrow
\infty}{1\over\Vol_d(rW)}{\Bbb E}\sum_{f\in\MaxFaces_k(Y(t\Lambda
,rW))}\varphi_k(f)\\
&=& \lim_{r\rightarrow\infty}{1\over\Vol_d(rW)}\int
\varphi_k(f){\Bbb F}_k^{Y(t\Lambda,rW)},
\end{eqnarray*}
not excluding thereby those maximal polytopes hitting the boundary of
$W\in\mathcal{K}_d$, which is somehow more natural in view of our
setting in Section \ref{secFIRST}.
\end{remark}

As a second example, we turn now to the length distribution of the
typical I-segment, which is nothing than the typical maximal polytope
of dimension one.
\begin{corollary}\label{corlength} The distribution of the length of
the typical I-segment of a stationary and isotropic STIT tessellation
with time parameter $t>0$ is a mixture of exponential distributions
with parameter $\gamma s$. The mixing distribution is a
beta-distribution on $(0,t)$ with parameters $d$ and $1$. Its density
is given by
\[
p_d(x) = \int_0^t\gamma s\mathrm{e}^{-\gamma sx}{ds^{d-1}\over t^d}\,\mathrm{d}s={d\over
(\gamma t)^dx^{d+1}}\Gamma(d+1,\gamma tx),\qquad x>0,
\]
where $\Gamma(\cdot,\cdot)$ is the lower incomplete Gamma-function
and $\gamma=\Gamma(d/2)/(\Gamma(1/2)\Gamma((d+1)/2))$.
\end{corollary}
\begin{pf} This follows immediately from Theorem \ref{CORDISTR} and
the well-known fact that the length distribution of the typical edge of
a stationary and isotropic Poisson hyperplane tessellation with
intensity $0<s<t$ is an exponential distribution with parameter $\gamma
s$, see~\cite{BL2}.
\end{pf}
In particular for $d=2$ and $d=3$, we have the densities
\begin{eqnarray*}
 p_2(x)
 &= &{1\over t^2x^3} \bigl(\uppi^2-(\uppi^2+2\uppi
tx+2t^2x^2)\mathrm{e}^{-(2/\uppi)tx} \bigr),\qquad x>0,\\
 p_3(x)
 &=& {3\over t^3x^4}
\bigl(48-(48+24tx+6t^2x^2+t^3x^3)\mathrm{e}^{-(1/2)tx}\bigr ),\qquad x>0.
\end{eqnarray*}
The mean segment lengths are $\uppi/t$ in the planar case and $3/t$ for
$d=3$. Moreover, the variance of the length of the typical I-segment in
the spatial case is given by $24/t^2$, which was not available before.
In general, from the explicit length density formula it is easily seen
that for the length of the typical I-segment only the moments of order
$1$ up to $d-1$ are finite.

Let us finally remark that Corollary \ref{corlength} allows an
extension to the anisotropic setting. Theorem \ref{CORDISTR} also
implies that the conditional length distribution of the typical
I-segment in $Y(t\Lambda)$ (where now $\Lambda$ is a general
translation-invariant hyperplane measure as in (\ref{LADEF})), given
its birth time $s\in(0,1)$ and direction $u\in{\Bbb S}^{d-1}$, is an
exponential distribution with parameter $s\Lambda([e(u)])$, where
$e(u)$ is a line segment of unit length parallel to $u$. Thus, the
conditional distribution of the length of the typical I-segment with a
given direction is a mixture of these exponential distributions and the
mixing distribution has again density ${ds^{d-1}\over t^d}\boldsymbol
{1}[0<s<t]$.
\begin{remark}
The length density $p_2(x)$ of the typical I-segment in a planar
stationary and isotropic STIT tessellation has been calculated in
\cite{MNW07} by an entirely different method based on Palm theory. However,
this method seems to be restricted to the study of I-segments and does
not lead to results for higher-dimensional maximal polytopes as in
Theorem \ref{CORDISTR} above.
\end{remark}

\section*{Acknowledgements}
We thank Joachim Ohser and Claudia Redenbach for providing the
simulations of the STIT tessellation shown in Figure \ref{Fig1}. The
second author would like to thank Werner Nagel and Matthias Reitzner
for their helpful comments and remarks. The suggestions and remarks of
two anonymous referees and of an associated editor were very helpful in
improving the presentation and the style of the paper.


\printhistory

\end{document}